\documentclass[12pt, a4paper]{article}
\usepackage[utf8]{inputenc}
\usepackage{amscd,amssymb,amsthm,amsmath}
\title{Goto's deformation theory of geometric structures, a Lie-theoretical description}
\author{Grigory Papayanov}
\date{2016}

\def\blacksquare{\hbox{\vrule width 7pt height 7pt depth 0pt}}
\def\endproof{\blacksquare}

\newcommand{\g}{{\mathfrak g}}
\newcommand{\arrow}{{\:\longrightarrow\:}}
\newcommand{\Z}{{\mathbb Z}}
\newcommand{\C}{{\mathbb C}}
\newcommand{\R}{{\mathbb R}}

\newcommand{\m}{\mathfrak m}
\def\1{\sqrt{-1}\:}
\renewcommand{\bar}{\overline}
\renewcommand{\phi}{\varphi}
\renewcommand{\epsilon}{\varepsilon}
\renewcommand{\geq}{\geqslant}

\newcommand{\End}{\operatorname{End}}
\newcommand{\Sym}{\operatorname{Sym}}
\newcommand{\Stab}{\operatorname{Stab}}
\newcommand{\Der}{\operatorname{Der}}
\newcommand{\Hom}{\operatorname{Hom}}
\newcommand{\Ann}{\operatorname{Ann}}

\newcommand{\Diff}{\operatorname{Diff}}

\newcommand{\Def}{\operatorname{Def}}
\newcommand{\Lie}{\operatorname{Lie}}

\newcommand{\Spec}{\operatorname{Spec}}
\newcommand{\MC}{\operatorname{MC}}
\newcommand{\Cone}{\operatorname{Cone}}

\newcommand{\A}{\mathcal{A}}
\newcommand{\EE}{\mathcal{E}}
\newcommand{\dR}{\Omega^\bullet}
\newcommand{\Ker}{\operatorname{Ker}}

\renewcommand{\exp}{\operatorname{exp}}

\renewcommand{\Im}{\operatorname{Im}}
\renewcommand{\Re}{\operatorname{Re}}

\newcounter{Mycounter}[section]
\newcounter{lemma}[section]
\renewcommand{\thelemma}{\noindent{Lemma \thesection.\arabic{lemma}}}
\newcommand{\lemma}{%
     \setcounter{lemma}{\value{Mycounter}}
     \refstepcounter{lemma}
     \stepcounter{Mycounter}
     {\bf \thelemma:\ }}

\newcounter{claim}[section]
\newcounter{sublemma}[section]
\newcounter{corollary}[section]
\newcounter{theorem}[section]
\renewcommand{\thetheorem}{\noindent{Theorem \thesection.\arabic{theorem}}}
\newcommand{\theorem}{%
     \setcounter{theorem}{\value{Mycounter}}
     \refstepcounter{theorem}
     \stepcounter{Mycounter}
     {\bf \thetheorem:\ }}	

\newcounter{conjecture}[section]
\newcounter{proposition}[section]
\newcounter{definition}[section]
\renewcommand{\thedefinition}
       {\noindent{Definition~\thesection.\arabic{definition}}}
\newcommand{\definition}{%
     \setcounter{definition}{\value{Mycounter}}
     \refstepcounter{definition}
     \stepcounter{Mycounter}
     {\bf \thedefinition:\ }}

\newcounter{example}[section]
\renewcommand{\theexample}{\noindent{Example \thesection.\arabic{example}}}
\newcommand{\example}{%
     \setcounter{example}{\value{Mycounter}}
     \refstepcounter{example}
     \stepcounter{Mycounter}
     {\bf \theexample:\ }}

\newcounter{remark}[section]
\renewcommand{\theremark}{\noindent{Remark \thesection.\arabic{remark}}}
\newcommand{\remark}{%
     \setcounter{remark}{\value{Mycounter}}
     \refstepcounter{remark}
     \stepcounter{Mycounter}
     {\bf \theremark:\ }}

\newcounter{problem}[section]
\newcounter{question}[section]
\begin{document}

\begin{center}
{\Large\bf
Goto's deformation theory of geometric structures, a Lie-theoretical description
}
\\[4mm]

Grigory Papayanov
\\[6mm]

\end{center}

%

\begin{abstract}

In \cite{Goto}, Ryushi Goto has constructed the deformation space
for a manifold equipped with a collection of closed differential forms
and showed that in some important cases (Calabi-Yau, $G_2$- and $Spin(7)$-structures)
this deformation space is smooth.
This result unifies the classical Bogomolov-Tian-Todorov and Joyce theorems
about unobstructedness of deformations.
 Using the work of Fiorenza and Manetti, we show that this deformation space
could be obtained as the deformation space associated to a certain $L_{\infty}$-algebra.
We also show that for Calabi-Yau, $G_2$- and $Spin(7)$-structures this
$L_{\infty}$-algebra is homotopy abelian. This gives a new proof of Goto's
theorem.

\end{abstract}




\section{Introduction}

The celebrated Berger's theorem (\cite{Besse}) classifies groups occuring as holonomy groups
of irreducible simply connedted Riemannian manifolds. The full list is

$$SO(n), U(n), Sp(n)\times Sp(1)/\Z_2, SU(n), Sp(n), G_2, Spin(7).$$

Manifold with four last holonomy groups are called, correspondingly,
Calabi-Yau manifolds, hyperK\"ahler manifolds, $G_2$-manifolds and $Spin(7)$-manifolds.
They share some interesting similar properties: for example, they are all Ricci-flat,
their cohomology admit Hodge decomposition into irreducible representations of the holonomy group,
and they have smooth local deformation spaces. For Calabi-Yau and hyperK\"ahler manifolds 
this smoothness was proved by Bogomolov (\cite{Bogomolov}), Tian (\cite{Tian}), and Todorov (\cite{Todorov}),
and for exceptional holonomy groups this was first proved by Joyce (\cite{Joyce}).
While Joyce's proof relies on complicated analytical arguments, 
proofs of Tian and Todorov are based on investigations of properties of dg-Lie algebra
$\Omega^{0,\bullet} \otimes T$ of tangent-valued forms.
The present work grew out of desire to have the similar Lie-algebraic description of
deformations of all manifolds with special holonomy.

\hfill

Other descriptions of moduli spaces of $G_2$ metrics belong to Hitchin (\cite{Hitchin})
and to Goto (\cite{Goto}). Hitchin's approach is based on finding extrema of a certain
functional defined on the space of forms, while Goto's approach is based on the following
observation. Suppose that
 tensor $A$ defines some type of geometric structure
(for example, in the sense of Cartan, i.e., the reduction of the structure group on the tangent
bundle from $GL$ to some smaller group stabilizing $A$), then
every other geometric structure (without any integrability conditions)
of this type is defined by $g\cdot A$ for some
$GL$-valued function on our base manifold. Thus, instead of deforming tensors, one
can deform the element $g$ of the gauge group.

\hfill

In (\cite{Goto}), Goto writes down conditions, under which it is possible to iteratively
solve the equation $d(g_t\cdot \alpha)=0$ for some closed form $\alpha$. Goto's arguments
are inductive and rather long. In the present work we are interpreting equations of (\cite{Goto})
as Maurer-Cartan equations in a certain $L_{\infty}$-algebra. The well-developed theory
of $L_{\infty}$-algebras allows to shorten Goto's calculations.

\hfill

{\bf Acknowledgements:} Author wants to thank D.Kaledin, M.Verbitsky and A.Zakharov for their 
interest in the present work and for the useful discussions. The study has been funded by
RSF grant, project 14-21-00053 dated 11.08.14.

%
%
%

\section{Geometric structures on manifolds}

In this section we define geometric structures and formulate the deformation problem.

\hfill

Let $V$ be a vector space, $\Phi \in \oplus_{n}\Lambda^{i_n}V^*$ be a collection of
exterior forms on it, and $\EE_{\Phi} := GL(V)\cdot \Phi$ be its orbit under the $GL(V)$-action
on the space $\oplus_{n}\Lambda^{i_n}V^*$. Then $\EE_{\Phi}$ is a homogeneous manifold isomorphic
to the $GL(V)/G$, where $G$ is the stabilizer of $\Phi$, with the tangent space at $\Phi$
isomorphic to the $\mathfrak{gl}(V)/\g$. 

\hfill

Now let $M$ be a differentiable manifold. Locally choosing for each point $x \in M$ an identification
between $T_xM$ and $V$, we can form the fiber bundle $\EE_{\Phi} \arrow M$,
consisting of forms that are pointwise in the $GL(V)$-orbit of $\Phi$. Denote by $\pi$
the natural projection from $\EE_{\Phi}$ to $M$.

\hfill

\definition
A $G$-{\bf structure}, or simply a geometric structure
(in the sense of Goto)
is the smooth section of the projection $\pi$. The geometric
structure is called {\bf integrable} if this section is a collection of closed forms.

\hfill

\remark
The notation is misleading in two ways. First, the subgroup $G$ of $GL(V)$
could be described by different antisymmetric
tensor invariants. Second, usually the name "$G$-structure" means the reduction
of the structure group of $TM$ from $GL(V)$ to $G$. The $G$-structure in
the sense of Goto provides the reduction of the structure group,
but not every reduction could be obtained by fixing skew-symmetric tensors.
To make things worse, Goto's notion of integrability does not always agree with other
possible notions.
However, as the following
examples show, this ambiguity should not present a problem, at least for the purposes of
this article.

\hfill

\example
An $Sp(2n)$-structure on an $2n$-dimensional manifold 
is the 2-form $\omega$ for which $\omega^n \ne 0$.
Integrable $Sp(2n)$-structure is called symplectic.

\hfill

\example
An $SL(n,\C)$-structure on an $2n$-dimensional manifold
is the complex $n$-form $\Omega$, for which
$\Omega \wedge \bar\Omega \ne 0$. An $SL(n,\C)$-structure on $M$ defines
an almost complex structure $I_{\Omega}$ on $M$, setting $T^{*1,0}M:=\{\alpha \in T^*M\mid\alpha\wedge\Omega=0\}$.
Integrability of the given $SL(n,\C)$ structure implies the integrability of
this complex structure.

\hfill

\example
An $SU(n,\C)$-structure on an $2n$-dimensional manifold
is the pair of forms $(\omega,\Omega) \in \Lambda^2M\oplus\Lambda^nM$,
where $\omega$ is an $Sp(2n)$-structure and $\Omega$ is the $SL(n,\C)$-structure,
$\Omega\wedge\omega=\bar\Omega\wedge\omega=0$,
$\Omega\wedge\bar\Omega=c\omega^n$ for some constant $c$, and the bilinear form
$\omega(I_{\Omega}\cdot,\cdot)$ is positive. A manifold
with an integrable $SU(n,\C)$-structure is called Calabi-Yau manifold.

\hfill

\example
Let $(\omega,\Omega)$ be an $SU(3)$-form on a 6-dimensional vector space $V$.
Consider the vector space $V \oplus \R$, and let $\theta$ be a coordinate on $\R$.
Consider forms
$$\phi=\omega\wedge\theta + \Im\Omega,~\psi=-\Re\Omega\wedge\theta + \frac{1}{2}\omega\wedge\omega.$$
The stabilizer of the pair $(\phi, \psi)$ is the exceptional Lie group $G_2$. 
For the 7-dimensional manifold $M$, sections of the corresponding
orbit $\EE_{(\phi,\psi)}$ are called $G_2$ structures. Integrable $G_2$ structures
are in 1-1 correspondence with Riemannian metrics of holonomy $G_2$ (\cite{Fernandez-Gray}).

\hfill

\remark
By the definition of $\EE_{\Phi}$, every section $\Phi_t$ could be obtained from
$\Phi$ by an action of the gauge group $GL(TM)$, $\Phi_t = g_t\cdot \Phi$,
with $g_t\cdot \Phi = h_t\cdot \Phi$ if and only if $g_th_t^{-1}$ lies in the stabilizer
of $\Phi$.

\hfill

%

\hfill

Since geometric structures are differential forms, pullbacks are defined; in particular,
the diffeomorphism group $\Diff(M)$ acts on the set of all integrable geometric structures.

\hfill

\definition
The set $\Gamma(\EE_{\Phi}) \cap \Ker(d) / \Diff(M)$ is called the {\bf moduli space} of $G$-structures.

\hfill

This is a factor of an infinite-dimensional manifold by an action of a 
Fr\'echet-Lie group, so at least it has the structure of a topological space. However, the neighbourhood of $\Phi$
admits a finer description. It makes sense to speak about infinitesimal neighbourhoods of $\Phi$
parametrized by local Artinian commutative $\R$-algebras.

\hfill

Let $\A$ be a local Artinian unital commutative
 $\R$-algebra with the maximal ideal $\m$. Denote by $\exp(\mathfrak{gl}(TM)\otimes \m)$ the
exponent of the nilpotent Lie algebra $\mathfrak{gl}(TM)\otimes \m$. 
The group $\exp(\mathfrak{gl}(TM)\otimes \m)$ acts on the 
sections of the bundle of $\A$-valued forms. If $\Phi$ is a section of $\Lambda^{\bullet}T^*M$,
then $\Phi \otimes \mathrm{1}$ is a section of $\Lambda^{\bullet}T^*M \otimes \A$, and we denote the
orbit of $\Phi \otimes \mathrm{1}$ under the action of $\exp(\mathfrak{gl}(TM)\otimes \m)$ by $\EE_{\Phi}(\A)$.
One can think of the sections of $\EE_{\Phi}(\A)$ as of the $\Spec(\A)$-parametrized families
of the sections of $\EE_{\Phi}$, with the central fiber isomorphic to the original form $\Phi$.
De Rham differential $d$ extends to $\Lambda^{\bullet}T^*(M) \otimes \A$ $\A$-linearly. We will
denote this extension also by $d$.

\hfill

Denote also by $\exp(\Gamma(TM) \otimes \m)$ the exponent of the nilpotent Lie algebra $\Gamma(TM) \otimes \m$.
The group $\exp(\Gamma(TM) \otimes \m)$ acts on the sections of $\EE_{\Phi}(\A)$ preserving closed sections.

\hfill

\definition\label{definition}
The functor from the category of local commutative Artinian $\R$-algebras to the category of sets,
defined by $$\Def_{\Phi}:(\A, \m) \mapsto \frac{(\EE_{\Phi} \otimes \A) \cap \Ker(d)}{\exp(\Gamma(TM)\otimes\m}$$ is 
called {\bf the deformation functor} of $\Phi$.

\hfill

\definition
Suppose $F$ is a functor from local commutative pro-Artinian $\R$-algebras to sets.
It is called {\bf unobstructed} if for every square-zero extension
$$0 \arrow I \arrow A' \arrow A \arrow 0$$

the induced map $F(A') \arrow F(A)$ is surjective.

\hfill

\theorem\label{Goto} (\cite{Goto})
A deformation functor associated to the Calabi-Yau, hyperK\"ahler, $G_2$- or $Spin(7)$-structure
is unobstructed.

\hfill

The unobstructedness of Calabi-Yau structures is a classical Bogomolov-Tian-Todorov theorem (\cite{Tian},\cite{Todorov}).
The unobstructedness of $G_2$- and $Spin(7)$-structures was proven by Joyce (\cite{Joyce}) with the help
of a lot of hard analysis. Our goal is to describe a deformation functor of a geometric structure as
a deformation functor associated to some dg-Lie algebra, giving a conceptual proof of Goto's theorem.

\hfill

\remark
Since one can exponentiate pro-nilpotent Lie algebras as well, the deformation functor of $\Phi$
is in fact defined on the category of local commutative pro-Artinian $\R$-algebras, 
such as the formal series algebra $\R[[t]]$. The question of convergence arises: is it true that the formal
$\Spec \R[[t]]$-parametrized deformation actually comes
 from the actual germ of a family of geometric structures?
If the manifold $M$ is compact, then the answer is often yes. 
The convergence of the solution of deformation equation considered in this article was proven by Goto in (\cite{Goto}).

\hfill

\section{DGLAs}

\definition
Let $(L,d,[\cdot,\cdot])$ be a differential graded Lie algebra with the 
inner derivation, that is,
there exists an element $\Delta \in L^1$ such that $dl=[\Delta,l]$ for every $l \in L$.
An element $x \in L^1$ is called {\bf Maurer-Cartan element} if $[\Delta+x,\Delta+x]=0$.
The set of all Maurer-Cartan elements is denoted by $\MC(L)$. 

\hfill

\definition
Suppose now that $L^0$ is a nilpotent Lie algebra. Denote its exponent by $\exp(L^0)$.
The adjoint action $ad$ of $L^0$ preserves the quadratic cone $\{l \in L^1 \mid [l,l]=0\}$,
and so is the action $e^{ad}$ of $\exp(L^0)$. The {\bf gauge action} is the action of
$\exp(L^0)$ on $\MC(L)$ given by $e^a * x := e^{ad(a)}(x+\Delta)-\Delta$.

\hfill

\remark
If the differential in $L$ is not inner, we can nevertheless force it to be so, by constructing
the new algebra $L'$. As a graded vector space,
$L'$ is equal to $L \oplus \langle \Delta \rangle$, where $\Delta$ is an element of degree 1.
Differential and commutators of $\Delta$ 
with other elements of $L$ are given by the relations
$d\Delta=0, [\Delta,l]=dl$. This construction
 allows one to define $\MC(L)$ and the gauge action for algebras with non-inner
derivation.

\hfill

\definition
The {\bf deformation functor} associated to $L$ is the functor $\Def_L$ from the category of local 
commutative Artinian rings to the category of sets given by

$$\Def_L(\A,\m):= \frac{\MC(L\otimes \m)}{\exp(L^0 \otimes \m)}.$$

\hfill

Here we list some well-known facts about deformation functors of dg-Lie algebras.
For the proofs, see e.g. (\cite{Manetti-lectures}).

\hfill

\theorem\label{ift}
Suppose $f: L \arrow M$ is a quasi-isomorphism of dg-Lie algebras. Then the induced morphism
 of functors $\Def(f): \Def_L \arrow \Def_M$ is an isomorphism.
 
\hfill

\definition
A dg-Lie algebra is called {\bf homotopy abelian} if it is quasi-isomorphic to the Lie algebra
with the zero bracket.

\hfill

\theorem\label{abelian}
If $M$ is  homotopy abelian, then the functor $\Def_M$ is unobstructed.
If $f:L \arrow M$ is a dg-Lie morphism, $\Def_M$ is unobstructed and $H^{\bullet}(f)$ is injective, then
$\Def_L$ is unobstructed as well.

%

\hfill

Next we describe the relative deformation functor construction, 
which was invented by Manetti (\cite{Manetti-relative})
and later explained by Manetti and Fiorenza (\cite{Manetti-Fiorenza}).

\hfill

Let $f:L\arrow M$ is the morphism of two dg-Lie algebras. Denote by $\MC(f)$ the set 
$$\{x \in L^1, e^m \in \exp(M^0)\mid x\in \MC(L), e^m*f(x)=0\}.$$
The group $\exp(L^0 \times M^{-1})$ acts on the set $\MC(f)$ by the rule $$(e^a,e^b)*(x,e^m)=(e^a*x,e^{db}e^me^{f(a)}).$$
This action is also called the gauge action.

\hfill		

\definition\label{reldeffun}
The {\bf relative deformation functor} associated to the morphism $f:L\arrow M$ is the functor from 
the category of local Artinian
algebras to the category of sets given by
$$\Def_f(\A,\m):=\frac{\MC(f\otimes \m)}{\mathrm{gauge~action}}.$$ 

\hfill
 
This functor was defined by Manetti in (\cite{Manetti-relative}) in order to
study deformations of a complex submanifold Lie-theoretically.
Later in (\cite{Manetti-Fiorenza}) Manetti and Fiorenza obtained
this functor as the deformation functor associated to a certain $L_{\infty}$-structure.
 We'll state their result, but for the sake of brevity we are not including 
the treatment of $L_{\infty}$-algebras into the present article. The reader 
can refer, for example, to the articles mentioned above.

\hfill

\definition
Let $f: V \arrow W$ be the map of complexes. Then {\bf the cone} of $f$ is the complex $\Cone(f)$
 with components $\Cone^i(f) := V^i \oplus W^{i-1}$, and the differential $d_f$, defined as
 $d_f(v,w):=(-dv,-f(v)+dw)$.
 
\hfill

\theorem(\cite{Manetti-Fiorenza})
Suppose $f: L \arrow M$ is the morphism of dg-Lie algebras. Define the operations 
$\langle \dots \rangle_n: \Sym^n(\Cone(f)) \arrow \Cone(f)$ by the following formulas:

\[
\langle (l,m)\rangle_1^{}=(-dl, -f(l)+dm);
\]
\[
\langle l_1\otimes l_2\rangle_2^{}=(-1)^{\deg_L(l_1)}[l_1,l_2];
\]
\[
\langle m\otimes
l\rangle_2^{}=\frac{(-1)^{\deg_M(m)+1}}{2}[m,f(l)];
\]
\[
\langle m_1\otimes
m_2\rangle_2^{}=0;
\]
\[
\langle{m_1}\otimes\cdots\otimes
{m_{n}}\otimes
l_1\otimes\cdots\otimes l_k\rangle_{n+k}^{}= 0, \qquad n+k\geq 3,
k\neq 1;
\]
and
\begin{multline*}
\langle {m_1}\otimes\cdots\otimes
{m_{n}}\otimes
l\,\rangle_{n+1}=\\
=-(-1)^{\sum_{i=1}^n\deg_M(m_i)}\frac{B_n}{n!}\sum_{\sigma\in
S_{n}}\varepsilon(\sigma)
[m_{\sigma(1)},[m_{\sigma(2)},\cdots,
[m_{\sigma(n)},f(l)]\cdots]],\qquad n\geq 2;
\end{multline*} 

where $B_n$ are Bernoulli numbers. Then these operations define 
the structure of an $L_{\infty}$-algebra on $\Cone(f)$, and
the functors $\Def_{\Cone(f)}$ and $\Def_f$ are isomorphic.
\endproof

\hfill

We will need the following analogues of \ref{ift} and \ref{abelian}.

\hfill

\theorem\label{properties} (\cite{Manetti-Fiorenza}, \cite{Manetti-relative})
Suppose $L,M,N,K$ are dg-Lie algebras, and the following square of dg-Lie algebra morphisms commutes: 

$$
\begin{CD}
L		@>i>>	N	\\
@VVfV		@VVgV	\\
M		@>j>>	K	\\
@VVV		@VVV	\\
\Cone(f)@.	 \Cone(g)
\end{CD}
$$

Then there exists a morphism $c: \Cone(f)\arrow \Cone(g)$, making the entire diagram commutative. Moreover,

\hfill

{\bf 1)} If $i$ and $j$ are quasi-isomorphims, then $c$ is also a quasi-isomorphism, and
functors $\Def_f$ and $\Def_g$ are isomorphic.

\hfill

{\bf 2)} If $\Cone(g)$ is homotopy abelian, then $\Def_g$ is unobstructed. Moreover,
if the map $H(c)$, induced by $c$ on cohomology groups, is an injection, then $\Def_f$ is also unobstructed.

%
%
%
%
%
%
%
%
%

%
%

\hfill

\remark
Suppose that $f: L\arrow M$ is an injection. Then for $(x,e^m)\in \MC(f)$, $x$ could be
reconstructed from $e^m$ by the formula $x=e^{-m}*0$. Hence, $\MC(f)$ for an inclusion could be
rewritten as $\{e^m \in \exp(M^0)\mid e^{-m}*0 \in L\}$.

\hfill


\hfill

\lemma\label{zerodiff}
Let $f: L \arrow M$ be the injection of graded Lie algebras 
(viewed as dg-Lie algebras with the trivial differential). Then $\Cone(f)$ is
homotopy abelian, and, consequently, $\Def_f$ is unobstructed.

\hfill

{\bf Proof:} 
It is proved in \cite{Martinengo-Fiorenza}, section 5.
 We present here another proof, based on the homotopy transfer theorem.
Reader interested in the homotopy transfer may refer, for example, to \cite{Loday-Valette}.

Let us pick a splitting $M = L \oplus M/L$, such that $f$ will be an 
embedding of the first summand, and denote the corresponding embedding of $M/L$ into $\Cone(f)$ by $e$.
Then the complex $M/L$ with the trivial differential is isomorphic to cohomology of the
complex $\Cone(f)$. Applying homotopy transfer formulas to this splitting, one can obtain
that transferred operations on $M/L$ vanish. Indeed, all brackets in $\Cone(f)$, applied to the elements
of $M \subset \Cone(f)$ vanish, and the image of $e$ lies in $M$. So, the transfered operations on the cohomology
of $\Cone(f)$ are zero, thus it is homotopy abelian. Hence the functor $\Def_{\Cone(f)}$ is unobstructed.
\endproof

\hfill

\theorem\label{matrix}
Let $(V,d)$ be a complex, $\End(V)$ its automorphism dg-Lie algebra
(that is, $\End^i(V):=\prod_{k \in \Z} \Hom(V^k,V^{k+i})$ and the differential is
given by the graded commutator with $d$), $v \in V$ 
 a closed (possibly non-homogeneous) vector, $dv=0$. Let $\Ann(v)$
be the subalgebra of endomorphisms annulating $v$, and denote the embedding 
$\Ann(v) \subset \End(V)$ by $\epsilon$. Then $\Def_{\epsilon}$ is unobstructed.

\hfill

{\bf Proof:} 
Let us pick a splitting of $V$ into an acyclic complex and a complex with the zero differential,
such that the projection of $v$ onto the acyclic summand is zero. This splitting induces splittings
of $\End(V)$ and $\Ann(v)$. Then the statement of the theorem follows
from \ref{properties} and \ref{zerodiff}.
\endproof

\section{Deformation of geometric structures}

Let $M$ be a smooth manifold, $\dR(M)$ its de Rham dg-algebra,
and $\Phi \in \Omega^{\bullet}(M)$ an integrable geometric structure on $M$. Our goal is
to describe functor from \ref{definition} as the deformation functor associated
to some $L_{\infty}$-algebra. We will describe it as the relative deformation functor.

\hfill

Denote by $\Der \dR(M)$ the dg-Lie algebra of derivations of $\dR(M)$. That is,
$\Der \dR(M)$ is the subcomplex of $\End \dR(M)$ consisting of those linear morphisms 
$D \in \prod_{k \in \Z} \Hom(\Omega^k, \Omega^{k+i})$ which satisfy the graded Leibniz rule: 
$D(\alpha \wedge \beta)=D(\alpha)\wedge \beta + (-1)^{k\cdot|\alpha|}\alpha \wedge D(\beta)$.

\hfill

\theorem\label{fronie}(\cite{Frolicher-Nijenhuis})
Every degree $k$ derivation on $\dR(M)$ is of the form $i_A+\{d,i_B\}$,
where $A \in \Omega^{k+1}(M,TM)$, $B \in \Omega^k(M,TM)$,
and $i_A$ is the inner derivation: $i_{\alpha\otimes v}(\beta):=\alpha \wedge i_v(\beta)$.
In particular, every derivation of degree $0$ is the sum of a
 Lie derivative with respect to some vector field $v$ and an infinitesimal gauge 
automorphism $A \in \mathfrak{gl}(TM)$. 

\hfill

Let $\Stab{\Phi}$ be a subalgebra of $\Der \dR(M)$ consisting of derivations vanishing on 
$\Phi$. Denote the inclusion $\Stab{\Phi} \arrow \Der \dR(M)$ by $j$.

\hfill

\theorem
The functor $\Def_j$ is isomorphic to the functor $\Def_{\Phi}$.

\hfill

{\bf Proof:} 
Remind that $MC_j(\m)$ is the set of those $m \in \Der(\Omega^0) \otimes \m$ for 
which $e^{\mathrm{ad}(-m)}d-d$ lies in $\Stab(\Phi)$, or, equivalently, $$e^{-m}de^m\Phi-d\Phi=0.$$
Since $\Phi$ is integrable, and $e^{-m}$ is invertible, this equation is equivalent to the
equation $de^m\Phi=0$. By \ref{fronie}, $m$ uniquely decomposes into a sum $A+\Lie_v$, where
$A \in \mathfrak{gl}(TM) \otimes \m$ and $v \in \Gamma(TM) \otimes \m$.

Consider the functor morphism $T: \MC_j \arrow \Def_{\Phi}$ given by the formula 
$$T(m)=T(A+\Lie_v)=e^A\Phi.$$ Since $A$ could be any element of $\mathfrak{gl}(TM)\otimes \m$,
this functor morphism is surjective. If $A$ and $B$ are two elements of 
$\mathfrak{gl}(TM)\otimes \m$, $e^A\Phi=e^B\Phi$ in $\Def_{\Phi}$ if and only if there exist
an element $v \in \Gamma(TM)\otimes\m$ such that $e^{-B}e^A\Phi=e^{Lie_v}\Phi$. This is
the equivalence under gauge action introduced in \ref{reldeffun}.
\endproof

\hfill

One advantage of working with $\Def_j$ rather than with $\Def_{\Phi}$ is the
access to \ref{ift} and \ref{abelian}.

Consider the following commutative diagram, where horizontal arrows are natural embeddings:

$$
\begin{CD}
\Stab(\Phi)	@>>>	\Ann(\Phi)	\\
@VVjV				@VViV		\\
\Der \dR(M)	@>>>	\End(\dR(M))
\end{CD}
$$

By naturality of Fiorenza-Manetti construction, this commutative square provides an $L_{\infty}$-map
from $\Cone(j)$ to $\Cone(i)$. By \ref{abelian}, $\Cone(i)$ is homotopy abelian, so the functor $\Def_i$
is unobstructed. If we will be able to prove that the induced map
$H^{\bullet}(\Cone(j)) \arrow H^{\bullet}(\Cone(i))$ is injective, we will prove the unobstructedness
of the functor $\Def_\Phi$.

\hfill

\theorem\label{main}
If $\Phi$ is a Calabi-Yau structure, hyperK\"ahler structure, $G_2$-structure or $Spin(7)$-structure,
then the map $H^{\bullet}(\Cone(j)) \arrow H^{\bullet}(\Cone(i))$ is injective.
\hfill

\hfill

{\bf Proof:} Both cones $\Cone(j)$ and $\Cone(i)$ are quasi-isomorphic to corresponding
quotients $\Der \dR(M)/\Stab(\Phi)$ and $\End(\dR(M))/\Ann(\Phi)$. Both factors are isomorphic
to orbits of $\Phi$ under the action of, correspondingly, $\Der \dR(M)$ and $\Ann(\Phi)$. Since the
action of $\End(\dR(M))$ can map any non-zero form onto any other non-zero form, 
the quotient 
$\End(\dR(M))/\Ann(\Phi)$ is isomorphic to the sum of several shifted de Rham complexes,
with the number of complexes being the number of non-zero components in $\Phi$
and the shifts being degrees of the corresponding components. Since $d\Phi=0$, from \ref{fronie} 
we obtain that the orbit of $\Phi$ under the action of $\Der \dR(M)$ is isomorphic to the following complex:

$$
0 \arrow i_{TM}\Phi \arrow \mathfrak{gl}(TM)\cdot\Phi \arrow \Omega^1\wedge\mathfrak{gl}(TM)\cdot\Phi \arrow 
\Omega^2\wedge\mathfrak{gl}(TM)\cdot\Phi \arrow \dots
$$

This is precisely the complex $\#$ of Goto (\cite{Goto}). Cohomological injectivity of the obvious map from
$\Der \dR(M)/\Stab(\Phi)$ to $\End(\dR(M))/\Ann(\Phi)$ in the cases of Calabi-Yau, hyperK\"ahler, $G_2$- and 
$Spin(7)$-structures are proven in chapters 4,5,6 and 7 of \cite{Goto}. In former two cases this injectivity 
is just a consequence of $\partial\bar\partial$-lemma, and in the latter, it is a consequence of some 
$\partial\bar\partial$-lemma-like statements about forms on manifolds with special holonomy.
\endproof

\hfill

%
%
%
%

%
%
%
%
%

\noindent {\sc {\bf G.P.}:
{\sc Laboratory of Algebraic Geometry,\\
National Research University Higher School of Economics,\\
Department of Mathematics, 6 Usacheva Str. Moscow, Russia,}\\
\tt  datel@mail.ru}.
\end{document}